\documentclass[12pt]{article}
\usepackage{hyperref}
\usepackage{amsfonts}
\usepackage{enumerate}
\usepackage{graphicx}
\usepackage{amsmath}
\begin{document}

\begin{center}
{\large\bf Causality of singular linear discrete time systems}

\vskip.20in
Nicholas Apostolopoulos$^{1}$, Fernando Ortega$^{2}$ and\ Grigoris Kalogeropoulos$^{3}$\\[2mm]
{\footnotesize
$^{1}$National Technical University of Athens, Greece\\
$^{2}$ Universitat Autonoma de Barcelona, Spain\\
$^{3}$National and Kapodistrian University of Athens, Greece}
\end{center}

{\footnotesize
\noindent
\textbf{Abstract:}
We consider two type of systems, a linear singular discrete time system and a linear singular fractional discrete time system whose coefficients are square constant matrices. By assuming that the input vector changes only at equally space sampling instants we investigate and provide properties for causality between state and inputs and causality between output and inputs.
\\
\\[3pt]
{\bf Keywords} : causality, linear, fractional, singular, system.
\\[3pt]

\vskip.2in

\section{Introduction}
Many authors have studied generalized discrete \& continouus time systems, see [1-28], and their applications, see [29-35]. Many of these results have already been extended to systems of differential \& difference equations with fractional operators, see [36-45]. 

If we define $\mathbb{N}_\alpha$ by $\mathbb{N}_\alpha=\left\{\alpha, \alpha+1, \alpha+2,...\right\}$, $\alpha$ integer, and $n$ such that $0<n<1$ or $1<n<2$, then the nabla fractional operator in the case of Riemann-Liouville fractional difference of $n$-th order for any 
$Y_k:\mathbb{N}_a\rightarrow \mathbb{R}^{m}$ is defined by
\[
\nabla_\alpha^{-n}Y_k=\frac{1}{\Gamma (n)}\sum^{k}_{j=\alpha}(k-j+1)^{\overline{n-1}}Y_j.
\]
We denote $\mathbb{R}^{m \times 1}$ with $\mathbb{R}^m$. Where the raising power function is defined by
\[
k^{\bar{\alpha}}=\frac{\Gamma(k+\alpha)}{\Gamma(k)}.
\]
We consider the singular discrete time system of the form
\begin{equation}
\begin{array}{cc}
FY_{k+1}=GY_k+V_k, & k= 1, 2,...,
\end{array}
\end{equation}
the singular fractional discrete time system of the form
\begin{equation}
\begin{array}{cc}
F\nabla_0^nY_k=GY_k+V_k, & k= 1, 2,...,
\end{array}
\end{equation}
the output system
\begin{equation}
\begin{array}{cc}
X_k=CY_k, & k= 1, 2,...,
\end{array}
\end{equation}
and the known initial conditions (IC)
\begin{equation}
Y_{0}.
\end{equation}      
Where  $F, G \in \mathbb{R}^{r \times m}$, $Y_k\in \mathbb{R}^{m}$, $V_k\in \mathbb{R}^{r}$, $X_k\in \mathbb{R}^{n}$ and $C \in \mathbb{R}^{n \times r}$. The matrices $F$, $G$ can be non-square ($r\neq m$) or square ($r=m$) with $F$ singular (det$F$=0). \\\\

\section{Preliminaries}
Throughout the paper we will use in several parts matrix pencil theory to establish our results. A matrix pencil is a family of matrices $sF-G$, parametrized by a complex number $s$, see [46-53].
\\\\
\textbf{Definition 2.1.} Given $F,G\in \mathbb{R}^{r \times m}$ and an arbitrary $s\in\mathbb{C}$, the matrix pencil $sF-G$ is called:
\begin{enumerate}
\item Regular when  $r=m$ and  det$(sF-G)\neq 0$;
\item Singular when  $r\neq m$ or  $r=m$ and det$(sF-G)\equiv 0$.
\end{enumerate}
In this article we consider the system (1) with a \textsl{regular pencil}, where the class of $sF-G$ is characterized by a uniquely defined element, known as the Weierstrass canonical form, see [46-53], specified by the complete set of invariants of $sF-G$. This is the set of elementary divisors of type  $(s-a_j)^{p_j}$, called \emph{finite elementary divisors}, where $a_j$ is a finite eigenvalue of algebraic multiplicity $p_j$ ($1\leq j \leq \nu$), and the set of elementary divisors of type $\hat{s}^q=\frac{1}{s^q}$, called \emph{infinite elementary divisors}, where $q$ is the algebraic multiplicity of the infinite eigenvalue. $\sum_{j =1}^\nu p_j  = p$ and $p+q=m$.
\\\\
From the regularity of $sF-G$, there exist non-singular matrices $P$, $Q$ $\in \mathbb{R}^{m \times m}$ such that 
\begin{equation}
\begin{array}{c}PFQ=\left[\begin{array}{cc} I_p&0_{p,q}\\0_{q,p}&H_q\end{array}\right],
\\\\
PGQ=\left[\begin{array}{cc} J_p&0_{p,q}\\0_{q,p}&I_q\end{array}\right].\end{array}
\end{equation}
$J_p$, $H_q$ are appropriate matrices with $H_q$ a nilpotent matrix with index $q_*$, $J_p$ a Jordan matrix and $p+q=m$. With $0_{q,p}$ we denote the zero matrix of $q\times p$. The matrix $Q$ can be written as
\begin{equation}
Q=\left[\begin{array}{cc}Q_p & Q_q\end{array}\right].
\end{equation}
$Q_p\in \mathbb{R}^{m \times p}$ and $Q_q\in \mathbb{R}^{m \times q}$. The matrix $P$ can be written as
\begin{equation}
P=\left[\begin{array}{c}P_1 \\ P_2\end{array}\right].
\end{equation}
$P_1\in \mathbb{R}^{p \times r}$ and $P_2\in \mathbb{R}^{q \times r}$.
The following results have been proved.
\\\\
\textbf{Theorem 2.1.}  (See [1-28]) We consider the systems (1) with a regular pencil. Then, its solution exists and for $k\geq 0$, is given by the formula
\[
    Y_k=Q_pJ_p^kD+QD_k.  
\]
Where $D_k=\left[
\begin{array}{c} \sum^{k-1}_{i=0}J_p^{k-i-1}P_1V_i\\-\sum^{q_{*}-1}_{i=0}H_q^iP_2V_{k+i}
\end{array}\right]$ and $D\in\mathbb{R}^p$ is a constant vector. The matrices $Q_p$, $Q_q$, $P_1$, $P_2$, $J_p$, $H_q$ are defined by (5), (6), (7). 
\\\\
\textbf{Theorem 2.2.} (See [36-45]) We consider the systems (2) with a regular pencil. Then, its solution exists if and only if all finite eigenvalues of the pencil are distinct and lie within the open disk $S=\{s\in\mathbb{R}:\left|s\right|<1\}$; Then, the solution of system (2) for $k\geq 0$, is given by the formula
\[
    Y_k=Q_p(k+1)^{\overline{n-1}}F_{n,n}(J_p(k+n)^{\bar n})(I_p-J_p)D+QD_k.  
\]
Where $D_k=\left[
\begin{array}{c} \sum^{k}_{i=1}(k-i+1)^{\overline{n-1}}F_{n,n}(J_p(k+n-i)^{\bar n})P_1V_i\\-\sum^{q_*-1}_{i=0}\nabla_0^{in}H_q^iP_2V_k
\end{array}\right]$ and $D\in\mathbb{R}^p$ is a constant vector. TThe matrices $Q_p$, $Q_q$, $P_1$, $P_2$, $J_p$, $H_q$ are defined by (5), (6), (7). The discrete Mittag-Leffler function with two parameters $F_{n,n}(J_p(k+n)^{\bar n})$ is defined by $F_{n,n}(J_p(k+n)^{\bar n})=\sum^{\infty}_{i=0}J_p^i\frac{(k+n)^{\overline{in}}}{\Gamma((i+1)n)}$.
\\\\
\textbf{Definition 2.2.} Consider the system (1) and (2)  with known IC of type (4). Then the IC are called consistent if there exists a solution for the system (1) and (2) respectively which satisfies the given conditions.
\\\\
\textbf{Proposition 2.1.} The IC of systems (1) and (2) are consistent if and only if 
\[
Y_0\in colspanQ_p +QD_0.
\]
\textbf{Proposition 2.2.}  Consider the system (1) with given IC. Then the solution for the initial value problem (1), (4) is unique if and only if the IC are consistent. Then, the unique solution is given by the formula
\[
    Y_k=Q_pJ_p^kZ^p_{0}+QD_k.  
\]
where $D_k=\left[
\begin{array}{c} \sum^{k-1}_{i=0}J_p^{k-i-1}P_1V_i\\-\sum^{q_{*}-1}_{i=0}H_q^iP_2V_{k+i}
\end{array}\right]$ and $Z^p_0$ is the unique solution of the algebraic system $Y_0=Q_pZ^p_0+D_0$.
\\\\
\textbf{Proposition 2.3.} Consider the system (1) with given IC (4). Then if there exists a solution for the initial value problem, it is unique if and only if the IC are consistent. Then, the unique solution is given by the formula
\[
    Y_k=Q_p(k+1)^{\overline{n-1}}F_{n,n}(J_p(k+n)^{\bar n})(I_p-J_p)Z^p_0+QD_k.  
\]
Where $D_k=\left[
\begin{array}{c} \sum^{k}_{i=1}(k-i+1)^{\overline{n-1}}F_{n,n}(J_p(k+n-i)^{\bar n})P_1V_i\\-\sum^{q_*-1}_{i=0}\nabla_0^{in}H_q^iP_2V_k
\end{array}\right]$ and $Z^p_0$ is the unique solution of the algebraic system $Y_0=Q_pZ^p_0+D_0$.

\section{Causality}

Generally for systems of type (1), (2) we define the notion of causality, which is properly defined bellow
\\\\
\textbf{Definition 3.1.} The singular systems (1), (2) are called causal, if its state $Y_k$, for any $k > 0$ is determined completely by initial state $Y_0$ and former inputs $V_0$, $V_1$, ..., $V_k$. Otherwise it is called non-causal.
\\\\
Next we will present the causality in the singular systems of the form (1), (2).
\\\\
\textbf{Proposition 3.1.} For system (1) causality between state and inputs exists if and only if there exists a matrix $B\in\mathbb{R}^{r\times r_1}$ such that $V_k=BU_k$ and $H_qP_2B=0_{q, r_1}$. Where $U_k\in\mathbb{R}^{r_1}$, while causality between output and input exists if and only if every column of the matrix $\left[\begin{array}{ccc}Q_qH_qP_2B&...&Q_qH_q^{q^*-1}P_2B\end{array}\right]$ lies in the right nullspace of the matrix $C$. The fractional singular systems of the form (2) is characterized by the property of causality. 
\\\\
\textbf{Proof.} From Proposition 2.2 the solution of system (1) with IC (4) is given by
\[
Y_k=Q_pJ_p^kZ^p_0+QD_k,
\]
or, equivalently,
\[
Y_k=Q_pJ_p^kZ^p_0+Q_p \sum^{k-1}_{i=0}J_p^{k-i-1}P_1V_i-Q_q\sum^{q_{*}-1}_{i=0}H_q^iP_2V_{k+i}.
\]
It is clear that causality of (1) depends on the term $\sum^{q_*-1}_{i=0}Q_qH_q^iP_2V_{k+i}$ and obviously $Y_k$ for any $k\geq 0$ is to be determined by former inputs if and only if there exists a matrix $B\in\mathbb{R}^{r\times r_1}$ such that $V_k=BU_k$ and $Q_qH_q^iP_2B=0_{q,r_1}$, or, equivalently, $H_q^iP_2B=0_{q,r_1}$, since $Q_p$ has linear independent columns it has a left inverse, or, equivalently, $H_qP_2B=0_{q, r_1}$. From (3) setting $Y_k$ in the state output $X_k$ we get
\[
X_k =CQ_pJ_p^{k-0}Z^p_{0}+CQ_p\sum^{k-1}_{i=0}J_p^{k-i-1}B_pV_i-CQ_q\sum^{q_{*}-1}_{i=0}H_q^iP_2V_{k+i}.
\]
From the above expression it is clear that non-causality is due to the existence of the term $\sum^{q_*-1}_{i=0}CQ_qH_q^iP_2V_{k+i}=\sum^{q_*-1}_{i=0}CQ_qH_q^iP_2BU_{k+i}$. So the causal relationship between $X_k$ and $V_k$ exists if and only if $CQ_qH_q^iP_2B=0_{m,r_1}$ for every $i= 1, 2, ..., q^*-1$. This can be written equivalently as $C\left[\begin{array}{ccc}Q_qH_qP_2B&...&Q_qH_q^{q^*-1}P_2B\end{array}\right]=0_{m,q^*nr_1}$ and thus system (1) is causal if and only if every column of the matrix $\left[\begin{array}{ccc}Q_qH_qP_2B&...&Q_qH_q^{q^*-1}P_2B\end{array}\right]$ lies in the right nullspace of the matrix $C$.
\\\\
From Proposition 2.3 the solution of system (2) with IC (4) is given by
\[
Y_k=Q_p(k+1)^{\overline{n-1}}F_{n,n}(J_p(k+n)^{\bar n})(I_p-J_p)Z^p_0+QD_k, 
\]
or, equivalently,
\[
Y_k=Q_pJ_p^kZ^p_0+Q_p\sum^{k}_{i=1}(k-i+1)^{\overline{n-1}}F_{n,n}(J_p(k+n-i)^{\bar n})P_1V_i-Q_q\sum^{q_*-1}_{i=0}\nabla_0^{in}H_q^iP_2V_k.
\]
It is clear that $Y_k$ of (2) is determined completely by the IC $Y_0$ ($Z^p_0$ is the unique solution of the linear system $Y_0=Q_pZ^p_0+D_0$) and the former inputs $V_0$, $V_1$, ..., $V_k$ and that occurs also between output and input of this system. The proof is completed. 

\section*{Conclusions}
In this article we focused and provided properties for causality of two types of system. A linear singular discrete time system and a linear singular fractional discrete time system whose coefficients are square constant matrices. 


\end{document}